\documentclass[11pt,letterpaper]{article}

\usepackage[utf8]{inputenc}
\usepackage{comment}
\usepackage{microtype}
\usepackage{graphicx}
\usepackage{appendix}
\usepackage{subcaption}
\usepackage{booktabs} 
\usepackage{scalerel}
\usepackage[margin=1in]{geometry}
\usepackage[most]{tcolorbox}
\usepackage{amssymb}
\usepackage{amsthm}
\usepackage{amsmath}
\usepackage{mathtools}
\usepackage{nicefrac}
\usepackage{algorithm}
\usepackage{algpseudocode} 
\usepackage{tabularx}
\usepackage[table]{xcolor}

\providecommand{ }{\hfill \ensuremath{\square}}
\usepackage{array}
\usepackage[hidelinks]{hyperref}

\usepackage[shortlabels]{enumitem}
\usepackage{dsfont}
\usepackage{url}
\usepackage{float}
\usepackage{pifont}
\usepackage{hhline}
\usepackage{multirow}
\usepackage{graphicx,wrapfig}
\usepackage{xcolor} 

\def\smskip{\smallskip}

\def\texitem#1{\par\smskip\noindent\hangindent 25pt
               \hbox to 25pt {\hss #1 ~}\ignorespaces}

\newcommand{\BEAS}{\begin{eqnarray*}}
\newcommand{\EEAS}{\end{eqnarray*}}
\newcommand{\BEA}{\begin{eqnarray}}
\newcommand{\EEA}{\end{eqnarray}}
\newcommand{\BEQ}{\begin{eqnarray}}
\newcommand{\EEQ}{\end{eqnarray}}
\newcommand{\BIT}{\begin{itemize}}
\newcommand{\EIT}{\end{itemize}}
\newcommand{\BNUM}{\begin{enumerate}}
\newcommand{\ENUM}{\end{enumerate}}

\newcommand{\BA}{\begin{array}}
\newcommand{\EA}{\end{array}}

\newif\ifpagenumbering
\pagenumberingtrue

\pagenumberingfalse

\newtheorem{theorem}{Theorem}

\newtheorem{corollary}[theorem]{Corollary}
\theoremstyle{remark}

\usepackage[numbers,square]{natbib}

\numberwithin{assumption}{section}
\numberwithin{definition}{section}
\numberwithin{theorem}{section}
\numberwithin{remark}{section}
\newcommand{\keywords}[1]{\textbf{Keywords:} #1}
\usepackage{hyperref}
\hypersetup{colorlinks,citecolor=blue,linktocpage,breaklinks=true}
\begin{document}
\title{Stability Analysis of Decentralized Adaptive Control of Laterally Coupled Multi-Loop Thermosyphon System with Unknown System Parameters}

\author{
\parbox{\textwidth}{
\centering
Novel Kumar Dey\thanks{Department of Applied Mathematics, The University of Arizona, Tucson, AZ, USA.}
\quad
Yan Wu\thanks{Department of Mathematical Sciences, Georgia Southern University, Statesboro, GA, USA.\\ 
Emails: \texttt{ndey@arizona.edu;yan@georgiasouthern.edu}}
}}

\date{} 
\maketitle
\begin{abstract}
{A high dimensional bi-directionally coupled $N$-loop thermosyphon system in a tandem is modeled by $3N$ Lorenz type of differential equations. The fluid flow in each loop is driven by the heat source (Rayleigh numbers) as well as the moment and heat exchanges between the adjacent loops. The flow becomes chaotic when the Rayleigh numbers are large. A decentralized controller design via proportional local state feedback is employed to stabilize the chaotic flows in each loop. As a result of the stability analysis, we show that there exist lower bounds on the controller gains that guarantee global stability of the system. Under the scenario of unknown parameters, we augment the original system with additional dynamic equations of the feedback gains to adaptively determine the feasible gains that stabilize the system. The analysis also sheds insight on the role of thermal coupling in the stability of the control system, which allows us to extend the results to negative $z$-state coupling coefficient as well as a class of nonlinear $z$-state coupling functions. Numerical simulations further demonstrate the effectiveness of the proposed controller design.}
\end{abstract}


\section{Introduction}\label{Intro}
Large-scale interconnected nonlinear dynamical systems arise in a broad range of modern engineering applications, including thermal-fluid transport networks, electric power systems, automated transportation, coordinated autonomous vehicles, manufacturing processes, and distributed sensing and actuation \cite{ref1,ref2}, and references therein. Such systems are composed of multiple interacting subsystems whose local dynamics are linked through physical, energetic, or informational exchanges. Although interconnection can enhance functionality and efficiency, it also creates substantial difficulties in analysis and control, particularly when the component dynamics are nonlinear, the system parameters are uncertain, and the number of subsystems is large. A central objective in nonlinear control is therefore to develop design and stability methodologies whose analytical structure remains tractable as the scale of the network increases \cite{ref1}.

\begin{figure}[H]
{\centering}{} 
\includegraphics[width=9.0cm]{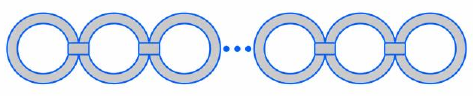}
    \caption{Laterally coupled multi-loop thermosyphon system.}\label{fig:laterally_coupled_n_loop}
\end{figure}

A paradigm of such a large-scale interconnected loop thermosyphon system and its control is studied in this paper. Multi-loop thermosyphon systems use gravity and natural convection across interconnected compartments to provide passive heat transfer. They are effective in managing high-density heat loads and offer multi-point thermal distribution. Their applications are found in nuclear decay heat removal systems, solar thermal arrays, and cooling systems for high-power servers stationed in large modern data centers. Consider a laterally coupled $N$-loop thermosyphon system, see Figure~\ref{fig:laterally_coupled_n_loop}, where the one-dimensional fluid flow in each loop is driven by the external heat sources imposed at the bottom of the loop and the heat sinks positioned at the top of the loop. Such a temperature difference drives the fluid in motion from regular convective flows to irregular time-dependent flows. There exist a heat flux and a flux of momentum exchanges between adjacent loops through the contact area, see Figure \ref{fig:laterally_coupled_n_loop}. Such exchanges allow one loop to interfere with the flow patterns in its neighboring loops, and vice versa.

There are three states associated with each loop, i.e., the state $x_i$ proportional to the average velocity in loop $i$ and the states $y_i$, $z_i$ proportional to the cross-sectionally averaged temperature distribution of the loop. The interconnected $N$-loop system is governed by the following system of $3N$ nonlinear differential equations,
\begin{align}\label{eq:general_n_loop_system}
    \begin{cases}
        \dot{x}_i
        =
        p\left\{
            y_i-x_i
            -g_{i-1}(x_i-x_{i-1})
            -g_i(x_i-x_{i+1})
        \right\},\\
        \dot{y}_i
        =
        R_i x_i-x_i z_i-y_i+u_i,\\
        \dot{z}_i
        =
        x_i y_i-z_i
        -h_{i-1}(z_i-z_{i-1})
        -h_i(z_i-z_{i+1}),
        \qquad i=1,2,\ldots,N.
    \end{cases}
\end{align}
where all system parameters in \eqref{eq:general_n_loop_system} are positive. The Rayleigh numbers, $R_i$, are the driven parameters. The system \eqref{eq:general_n_loop_system} becomes chaotic under large Rayleigh numbers. The momentum coupling intensity coefficient between loop $i$ and the loop after is given by function $g_i$ with the convention $g_0=g_N=0$. Similarly, $h_i$ is the thermal coupling intensity function between the corresponding adjacent loops with $h_0=h_N=0$. Lastly, the term $u_i$ represents the control signal. Our control objective is to stabilize \eqref{eq:general_n_loop_system} well into its chaotic regime under the assumption that all system parameters are unknown and an arbitrary number of loops being considered, particularly, when the number of loops $N$ is large. It is noted that the interconnected, large scale system \eqref{eq:general_n_loop_system} is the result of a generalization of the two-loop system first reported in \cite{ref3}, where detailed bifurcation analysis of the periodic solutions can be found. Previous work on coupled loop thermosyphons has treated several fixed low-dimensional configurations. Decentralized proportional state feedback control has been used to stabilize two-loop \cite{ref4}, three-loop \cite{ref5}, and four-loop systems \cite{ref6,ref7}, and adaptive gain laws have been introduced when the Rayleigh numbers and coupling parameters are unknown. The recent triangular three-loop study \cite{ref8} further incorporated extended-state observers and active disturbance rejection control and established explicit gain bounds that depend on the Rayleigh numbers and momentum coupling parameters but not on the positive thermal coupling coefficients. These studies provide important evidence that local temperature feedback is both analytically effective and practically attractive. Nevertheless, their stability proofs rely on topology-specific examination of low-dimensional Lyapunov matrices and direct expansion of their principal minors. Furthermore, a common denominator of the controller design in those papers is decentralized control, where the controller for each loop only depends on the local state, i.e., $u_i=-k_i y_i$, $i=1,2,\ldots,N$, where $k_i$'s are the positive feedback gains.

In the domain of interconnected control systems, decentralized control is especially well suited to this setting because each local controller uses only locally available measurements and does not require access to the complete network states. This architecture reduces communication and implementation burdens, improves modularity and fault tolerance, and allows additional subsystems to be incorporated without redesigning a centralized controller. Recent investigations have extended decentralized adaptive control to interconnected systems with non-triangular structural uncertainty and time-varying parameters \cite{ref9}, prescribed transient performance \cite{ref10}, unknown disturbances and finite-time stabilization objectives \cite{ref11}, and time delays, actuator failures, and event-triggered output-feedback \cite{ref12}. Related event-triggered designs have also addressed input saturation and finite-time convergence in uncertain interconnected nonlinear systems \cite{ref13}. These developments demonstrate the breadth of contemporary decentralized control theory, while also illustrating that much of the recent literature has concentrated on increasingly sophisticated controller architectures for broad abstract system classes.

An equally important question is scalability of the accompanying stability analysis. For a fixed interconnected system, positive definiteness of a Lyapunov matrix may be checked numerically, and a sufficiently large feedback gain may often be selected by trial and error. Such procedures, however, provide limited analytical insight and generally do not yield explicit controller-design conditions that remain interpretable as the system dimension grows. In physically structured networks, the central challenge is not only to prove stability, but also to exploit sparsity and interconnection topology so that the stability conditions can be generated systematically for an arbitrary number of compartments. This issue is particularly significant when the local dynamics are nonlinear, and the coupling parameters have direct physical meaning. Interconnected thermal-fluid systems provide an important setting in which these questions arise. Recent studies of multi-producer, district-heating networks, for example, have developed decentralized adaptive and passivity-based controllers for regulating flows, storage volumes, and network temperatures while preserving closed-loop stability \cite{ref14,ref15}. These results confirm that modular control strategies can be effective for high-dimensional thermal networks with distributed components. The thermosyphon problem considered here is related in some respects but differs substantially in its local dynamics: each loop may undergo a transition from conduction to steady convection, periodic oscillation, and chaotic flow as the Rayleigh number increases. When several loops exchange momentum and heat, the resulting coupled system combines strongly nonlinear Lorenz-type dynamics with a sparse nearest-neighbor interconnection structure.

Loop thermosyphons transfer heat through buoyancy-driven circulation without mechanical pumps and have long served as representative nonlinear thermal-fluid systems \cite{ref3}. In a single loop, the Rayleigh number measures the strength of the imposed heating and plays a decisive role in the onset of oscillatory and chaotic behavior. In a coupled configuration, momentum exchange alters the fluid velocities, whereas thermal exchange links the temperature states of adjacent loops. The interaction of these mechanisms produces a high-dimensional dynamical system in which chaotic motion in one compartment can influence the neighboring compartments. From a control perspective, the problem is therefore to suppress the chaotic flows using a controller architecture that is physically measurable, decentralized, and scalable with the number of loops.

This dependence on a fixed topology as seen in \cite{ref4,ref5,ref6,ref7,ref8} becomes an obstacle when the number of loops is arbitrary. A tandem network of $N$ thermosyphon loops is governed by $3N$ nonlinear differential equations, and each additional loop enlarges the Lyapunov matrix as seen in \cite{ref4,ref5,ref6,ref7,ref8} while introducing new coupling terms. Although Sylvester's criterion gives a complete characterization of positive definiteness, direct determinant calculations rapidly lose transparency and become extremely cumbersome. The formulas obtained for one network size do not automatically reveal the corresponding conditions for the next size, and the physical role of the thermal and momentum couplings becomes increasingly obscured. A general theory therefore requires a structural reduction that separates the different coupling mechanisms and a recursive procedure that preserves its form as the number of loops increases. The present paper develops such a recursive principal-minor stability framework for a laterally coupled $N$-loop thermosyphon network. The system is modeled by a $3N$-dimensional Lorenz-type set of equations and is controlled through decentralized proportional feedback based on the local horizontal-temperature state. By applying an appropriate permutation to the state variables, the resulting symmetric Lyapunov matrix is transformed into a block-diagonal form consisting of a thermal block and a reduced momentum-temperature block. The thermal block is shown to be positive definite under the standard nonnegative coupling intensity assumption. As such, the controller-gain conditions are determined entirely by the reduced sparse block-tridiagonal matrix.

The fixed-gain stability theory is then used as the foundation for decentralized adaptive control under parameter uncertainty. Because the explicit thresholds cannot be evaluated when the Rayleigh numbers or coupling parameters are unavailable, each feedback gain is augmented with a local dynamic law that increases the gain while the corresponding temperature state remains nonzero. The existence of finite stabilizing thresholds guarantees that the adaptive gains eventually enter a stabilizing region, after which the physical states converge to the equilibrium. The analysis further reveals that the thermal block can remain positive definite beyond the usual positive linear-coupling setting. This observation allows the results to be extended to negative coupling within an admissible range and to a class of nonlinear thermal coupling functions.

The main contributions of this work may be summarized as follows. First, a unified $3N$-dimensional model and stability analysis framework are established for an arbitrary tandem of laterally coupled thermosyphon loops. Second, explicit conditions on the feedback gains that guarantee global stability of the $N$-loop system are derived, which allows one to construct adaptive gain laws for unknown system parameters. Third, a permutation-based matrix decomposition and recursive principal-minor method convert a dimension-dependent positive-definiteness problem into a systematic procedure that scales naturally with the network size. The recursive relations of the principal minors may be utilized to compute less conservative bounds on the gains when the system parameters are known.  Finally, the role of thermal coupling is characterized in a form that supports extensions to negative and nonlinear coupling mechanisms. Numerical simulations are provided to underpin the theoretical results and showcase the performance of the proposed decentralized adaptive controller design. The rest of the paper is organized as follows: in Section \ref{sec:Coupled_$N$_Loop_System_And_Adaptive_Control}, we introduce a unique decoupled structure of the stability matrix associated with the $N$-loop system. We then present the main results on the stability of the $N$-loop system in two cases: the existence of lower bounds on the feedback gains and the adaptive control with unknown system parameters. In the same section, we derive recursive relations between the principal minors of the stability matrix, which allows step-by-step computation of the lower bounds on the decentralized controller gains. We discuss nonlinear coupling in Section \ref{sec:nonlinear_thermal_coupling}. The unique structure of the permuted stability matrix allows us to extend the linear thermal coupling to a more general setting such as negative coupling or nonlinear coupling. We derive explicit conditions on the nonlinear coupling function, which retain the stability of the multi-loop system. The main results and conclusive remarks are summarized in Section \ref{sec:conclusions}.

\section{Generalized Laterally Coupled Multi-Loop System}\label{sec:Coupled_$N$_Loop_System_And_Adaptive_Control}
We consider the generalized laterally coupled $N$-loop thermosyphon system, $N\geq2$, with constant momentum and thermal coupling parameters. We then rewrite \eqref{eq:general_n_loop_system} into the following form:
\begin{align}\label{eq:n_sys}
    \begin{cases}
    \dot{x}_i
    =
    p\{(y_i-x_i)-\gamma_{i-1}(x_i-x_{i-1})-\gamma_i(x_i-x_{i+1})\} \\
    \dot{y}_i
    =
    R_i x_i-x_i z_i-k_i y_i \\
    \dot{z}_i
    =
    x_i y_i-z_i-\eta_{i-1}(z_i-z_{i-1})-\eta_i(z_i-z_{i+1}),
    \qquad i=1,\ldots,N,
    \end{cases}
\end{align}
with the decentralized single-state feedback controller $u_i=-k_i y_i,
\quad i=1,\ldots,N$. We identify the feedback gains $k_i$ as the driving parameters for stabilizing the system \eqref{eq:n_sys}. For the purpose of stability analysis, we assemble the state vector from the outset by placing all thermal $z$-states at the end:
\begin{align}\label{eq:ordered_state_vector}
    X=(x_1,y_1,x_2,y_2,\ldots,x_N,y_N, z_1,z_2,\ldots,z_N)^\top.
\end{align}
Thus, the first $2N$ components of $X$ contain the $(x,y)$-states, while the final $N$ components contain the $z$-states. The endpoint conventions are $ \gamma_0=\gamma_N=0,
    \quad \eta_0=\eta_N=0.$
In the uniform laterally coupled case considered here, we take $\gamma_i=\gamma,
    \quad
    \eta_i=\eta,
    \quad
    i=1,\ldots,(N-1).$ Throughout this section, we assume
\begin{align}\label{eq:standing_parameter_assumptions}
    R_i>0,
    \qquad
    0<\gamma<1,
    \qquad
    \eta\geq0,
    \qquad
    i=1,\ldots,N.
\end{align}
The endpoint conventions ensure that all nonexistent boundary terms vanish at $i=1$ and $i=N$. To analyze the stability of \eqref{eq:n_sys}, we introduce the symmetric stability matrix
$A_N(k)\in\mathbb{R}^{3N\times3N}$.
With the state ordering in \eqref{eq:ordered_state_vector},
this matrix has the block-diagonal representation of the following form
\begin{align}\label{eq:AN_block_diagonal_form}
    A_N(k)
    =
    \begin{bmatrix}
        H_N(k) & 0\\
        0 & T_N(\eta)
    \end{bmatrix}.
\end{align}
Here, $H_N(k)\in\mathbb{R}^{2N\times2N}$ contains the contribution of the $(x,y)$-states, and $T_N(\eta)\in\mathbb{R}^{N\times N}$ contains the contribution of the thermal $z$-states. The reduced block $H_N(k)$ is ordered according to $(x_1,y_1,x_2,y_2,\ldots,x_N,y_N)$ and has the block-tridiagonal form
\begin{align}\label{eq:HN_block_tridiagonal_form}
    H_N(k)
    =
    \begin{bmatrix}
        B_1       & C_1       & 0          & \cdots & 0\\
        C_1^\top  & B_2       & C_2        & \ddots & \vdots\\
        0         & C_2^\top  & B_3        & \ddots & 0\\
        \vdots    & \ddots    & \ddots     & \ddots & C_{N-1}\\
        0         & \cdots    & 0          & C_{N-1}^\top & B_N
    \end{bmatrix},
\end{align}
where the block matrix associated with the $i$-th loop in \eqref{eq:n_sys} is given by
\begin{align}\label{eq:local_Bi_block}
    B_i
    =
    \begin{bmatrix}
        d_i & -R_i\\
        -R_i & k_i
    \end{bmatrix},
    \qquad i=1,\ldots,N,
\end{align}
where
\begin{align}\label{eq:di_definition}
    d_i
    =
    \begin{cases}
        R_i(1+\gamma),
        & i=1 \text{ or } i=N,\\[1mm]
        R_i(1+2\gamma),
        & 2\leq i\leq (N-1).
    \end{cases}
\end{align}
The coupling block between adjacent loops $i$ and $i+1$ is
\begin{align}\label{eq:Ci_block_definition}
    C_i
    =
    \begin{bmatrix}
        -q_i & 0\\
        0 & 0
    \end{bmatrix},
    \qquad
    q_i
    = (R_i+R_{i+1})\gamma/2,
    \qquad i=1,\ldots,(N-1).
\end{align}
Therefore, the only inter-loop entries in $H_N(k)$ are the couplings $-q_i$ between the neighboring velocity states $x_i$ and $x_{i+1}$. Each $y_i$-state couples only with its corresponding $x_i$-state through the entry $-R_i$.

The thermal block is
\begin{align}\label{eq:thermal_block_definition}
    T_N(\eta)
    =
    \begin{bmatrix}
        1+\eta & -\eta & 0 & \cdots & 0\\
        -\eta & 1+2\eta & -\eta & \ddots & \vdots\\
        0 & -\eta & 1+2\eta & \ddots & 0\\
        \vdots & \ddots & \ddots & \ddots & -\eta\\
        0 & \cdots & 0 & -\eta & 1+\eta
    \end{bmatrix}.
\end{align}

Thus, the stability conditions on the feedback gains $k_i$ are determined entirely by the
reduced blocks $H_N(k)$ and $T_N(\eta)$.

\subsection{Global Asymptotic Stability of the $N$-Loop System}
Since matrix $A_N(k)$ is in a block-diagonal form, it is positive definite if and only if the submatrices $H_N(k)$ and $T_N(\eta)$ are positive definite. The following theorem establishes the existence of lower bound on the feedback gains that guarantee $A_N(k)$ is positive definite, which proves the global asymptotic stability of the closed-loop system \eqref{eq:n_sys}.

\begin{theorem}\label{thm:positive_definiteness_before_recursions}
Consider the laterally coupled $N$-loop system
\eqref{eq:n_sys} with constant linear coupling. Suppose that
\begin{align}\label{eq:theorem_rayleigh_closeness}
    \gamma
    \left(
        \frac{\max_{1\leq j\leq N}R_j}
             {\min_{1\leq j\leq N}R_j}
        -1
    \right)
    <1.
\end{align}
There exists lower bound on the feedback gains $k_i$ such that the close loop system \eqref{eq:n_sys} is globally asymptotically stable at the origin.
\begin{proof}
Consider the Lyapunov function
\begin{align}\label{eq:lyapunov_function_general_N}
    V(X)
    =
    \frac12\sum_{i=1}^{N}
    \left(
        \frac{R_i}{p}x_i^2+y_i^2+z_i^2
    \right).
\end{align}
Since $p>0$ and $R_i>0$, this function is continuously differentiable, positive, and radially unbounded. Differentiating $V$ along the trajectories of
\eqref{eq:n_sys} gives
\begin{align*}
    \dot V(X)
    =
    \sum_{i=1}^{N}
    \left(
        \frac{R_i}{p}x_i\dot x_i
        +y_i\dot y_i
        +z_i\dot z_i
    \right),
\end{align*}
which can be written as a quadratic form

\begin{align}\label{eq:lyapunov_derivative_matrix_form}
    \dot V(X)
    &=
    -\sum_{i=1}^{N}d_i x_i^2
    +2\sum_{i=1}^{N}R_i x_i y_i
    -\sum_{i=1}^{N}k_i y_i^2+
    2\sum_{i=1}^{N-1}q_i x_i x_{i+1} \notag\\
    & \qquad-\sum_{i=1}^{N}z_i^2
    -\eta\sum_{i=1}^{N-1}(z_i-z_{i+1})^2
    \notag\\
    &=
    -X^\top A_N(k)X,
\end{align}
where the state vector $X= (x_1,y_1,\ldots,x_N,y_N,z_1,\ldots,z_N)^\top$. The matrix  $A_N(k)$ is given by \eqref{eq:AN_block_diagonal_form} as well as it's submatrices $H_N(k)$ \eqref{eq:HN_block_tridiagonal_form} and $T_N(\eta)$ \eqref{eq:thermal_block_definition}. In what follows, we will derive conditions for $H_N(k)$ and $T_N(\eta)$ to be positive definite. For every interior index $i=2,\ldots,N-1$, using \eqref{eq:di_definition} and \eqref{eq:Ci_block_definition}, we obtain
\begin{align*}
    d_i-q_{i-1}-q_i
    &=
    R_i(1+2\gamma)
    -
    \frac{(R_{i-1}+R_i)\gamma}{2}
    -
    \frac{(R_i+R_{i+1})\gamma}{2}\\
    &=
    R_i(1+\gamma)
    -
    \frac{\gamma}{2}
    \left(
        R_{i-1}+R_{i+1}
    \right).
\end{align*}
Because $R_{i-1},R_{i+1}
    \leq
    \frac{
        \displaystyle\max_{1\leq j\leq N}R_j
    }{
        \displaystyle\min_{1\leq j\leq N}R_j
    }
    R_i,$ it follows that
\begin{align*}
    d_i-q_{i-1}-q_i
    &\geq
    R_i
    \left[
        1-
        \gamma
        \left(
            \frac{
                \displaystyle\max_{1\leq j\leq N}R_j
            }{
                \displaystyle\min_{1\leq j\leq N}R_j
            }
            -1
        \right)
    \right]
    >0,
\end{align*}
by \eqref{eq:theorem_rayleigh_closeness}. Corresponding to the first loop, we have
\begin{align*}
    d_1-q_1
    &=
    R_1(1+\gamma)
    -
    \frac{(R_1+R_2)\gamma}{2}\\
    &=
    R_1
    +
    \frac{\gamma}{2}(R_1-R_2)\\
    &\geq
    R_1
    \left[
        1-
        \frac{\gamma}{2}
        \left(
            \frac{
                \displaystyle\max_{1\leq j\leq N}R_j
            }{
                \displaystyle\min_{1\leq j\leq N}R_j
            }
            -1
        \right)
    \right]
    >0.
\end{align*}
The same argument applied to the final loop gives $d_N-q_{N-1}>0.$ Therefore,
\begin{align}\label{eq:theorem_positive_structural_terms}
    d_i-q_{i-1}-q_i>0,
    \qquad i=1,\ldots,N,
\end{align}
with the notion $q_o=q_N=0$. 
We will prove that $H_N(k)$ is positive definite because of \eqref{eq:theorem_positive_structural_terms}. Let $\xi
    =
    (x_1,y_1,x_2,y_2,\ldots,x_N,y_N)^\top
    \in\mathbb{R}^{2N}.$ From the structure of $H_N(k)$, we obtain
\begin{align*}
    \xi^\top H_N(k)\xi
    =
    \sum_{i=1}^{N}d_i x_i^2
    -
    2\sum_{i=1}^{N-1}q_i x_i x_{i+1}-
    2\sum_{i=1}^{N}R_i x_i y_i
    +
    \sum_{i=1}^{N}k_i y_i^2.
\end{align*}
Completing the squares in $y_i$ gives
\begin{align*}
    \xi^\top H_N(k)\xi
    =
    \sum_{i=1}^{N}
    k_i
    \left(
        y_i-\frac{R_i}{k_i}x_i
    \right)^2+
    \sum_{i=1}^{N}
    \left(
        d_i-\frac{R_i^2}{k_i}
    \right)x_i^2
    -
    2\sum_{i=1}^{N-1}q_i x_i x_{i+1}.
\end{align*}
Since $-2q_i x_i x_{i+1}
    \geq -q_i x_i^2-q_i x_{i+1}^2,$ we obtain
\begin{align*}
    \xi^\top H_N(k)\xi
    \geq
    \sum_{i=1}^{N}
    k_i
    \left(
        y_i-\frac{R_i}{k_i}x_i
    \right)^2+
    \sum_{i=1}^{N}
    \left(
        d_i-q_{i-1}-q_i-\frac{R_i^2}{k_i}
    \right)x_i^2.
\end{align*}
By setting $d_i-q_{i-1}-q_i-\frac{R_i^2}{k_i}>0$, and by \eqref{eq:theorem_positive_structural_terms}, we have
\begin{align}\label{eq:theorem_preliminary_gain_bounds}
    k_i
    >
    \frac{R_i^2}{d_i-q_{i-1}-q_i},
    \qquad i=1,\ldots,N.
\end{align}
Consequently, $\xi^\top H_N(k)\xi>0$ for every nonzero $\xi$. Hence,
\begin{align}\label{eq:theorem_HN_positive_definite}
    H_N(k)\succ0.
\end{align}
As far as the thermal block $T_N(\eta)$ is concerned, for every nonzero vector $z=(z_1,\ldots,z_N)^\top\in\mathbb{R}^N,$ we have
\begin{align*}
    z^\top T_N(\eta)z
    =
    \sum_{i=1}^{N}z_i^2
    +
    \eta\sum_{i=1}^{N-1}(z_i-z_{i+1})^2.
\end{align*}
Since $\eta>0$, it follows that
\begin{align}\label{eq:theorem_TN_positive_definite}
    T_N(\eta)\succ0.
\end{align}
Hence, $A_N(k)
    = \begin{bmatrix}
        H_N(k) & 0\\
        0 & T_N(\eta)
    \end{bmatrix}$ is positive definite, which gurantees $
    \dot V(X)<0
    \quad\text{for all }X\neq0.$ Therefore, there exist lower bounds on $k_i$ satisfying \eqref{eq:theorem_preliminary_gain_bounds} such that the $N$-loop system \eqref{eq:n_sys} is globally asymptotically stable at the origin. 
\end{proof}
\end{theorem}

When the system parameters are unknown, we extend the state space with the $y$-state
feedback gains as new state variables along with their own dynamic equations so that
they can be determined adaptively.  The augmented system of differential equations is
given by

\begin{align}\label{eq:aug_system}
\begin{cases}
\dot{x}_i &= p\{y_i-x_i-\gamma_{i-1}(x_i-x_{i-1})
-\gamma_i(x_i-x_{i+1})\} \\
\dot{y}_i &= R_i x_i-x_i z_i-k_i y_i \\
\dot{z}_i &= x_i y_i-z_i-\eta_{i-1}(z_i-z_{i-1})
-\eta_i(z_i-z_{i+1}) \\
\dot{k}_i &= \beta_i y_i^2  \qquad i=1,2,\ldots,N.
\end{cases}
\end{align}
where $\beta_i>0$ and $k_i(0)=0$, $i=1,2,\ldots,N$.  It is readily seen from \eqref{eq:aug_system} that the extended
state $k_i(t)$ is monotonically increasing in time. As a result of Theorem \ref{thm:positive_definiteness_before_recursions}, each
gain $k_i$ eventually reaches their respective threshold value $k_i^*$, so that system \eqref{eq:aug_system} is globally
asymptotically stable.  Therefore, the equilibrium of \eqref{eq:aug_system} is identified as
$(\underline{0},k^*)$, where $k^*=(k_i^*)_{i=1}^N$.  The stability result of system \eqref{eq:aug_system} is
summarized in the following theorem.

\begin{theorem}\label{thm:global}
The extended closed-loop system \eqref{eq:aug_system} is globally asymptotically stable at
the equilibrium.
\begin{proof}Consider the Lyapunov function,

\begin{equation}\label{eq:aug_system_Lyap}
V=\frac{1}{2}\sum_{i=1}^{N}
\left(
\frac{R_i}{p}x_i^2+y_i^2+z_i^2+
\frac{1}{\beta_i}(k_i-k_i^*)^2
\right)
\end{equation}

The time derivative $\dot{V}$ along the state trajectory is written as
$\dot{V}=-w^T A(k^*)w$, where $w$ is the same state vector as in the proof of Theorem \ref{thm:positive_definiteness_before_recursions}, and $A(k^*)$ is the matrix given by \eqref{eq:AN_block_diagonal_form}, where $k_i$ is replaced by
$k_i^*$.  Since the gain states $k_i$'s are absent from $\dot{V}$, we can only
conclude $\dot{V}\leq 0$.  We apply Barb\u{a}lat's lemma to prove that
$w\to\underline{0}$ as follows. First, we show $w$ is uniformly continuous.  To this
end, it is sufficient to show that $\dot{w}$ is bounded. It is known that, in general, $R_i>p$ in \eqref{eq:aug_system_Lyap}, which implies
\begin{align*}
    V&>\frac{1}{2}\sum_{i=1}^{N}x_i^2+y_i^2+z_i^2 =\|w\|^2/2.
\end{align*}    
Since
$\dot{V}\leq 0$, one has $V(t)\leq V(0)$.  Hence, $\|w\|^2\leq 2V(0)$.
Therefore, all states $x_i$, $y_i$, and $z_i$ are bounded.  

Similarly, it can be
shown that $\sum_{i=1}^{N}(k_i-k_i^*)^2\leq 2\beta V(0)$, where
$\beta=\max_i\{\beta_i\}$. This shows that each $k_i$ is bounded. We then conclude
that $\dot{w}$ is bounded based on equations in \eqref{eq:aug_system}. Therefore, $w$ is uniformly continuous. 

Next, we show that $\int_0^\infty\|w\|^2\,dt$ exists and finite. Again,
from $\dot{V}=-w^T A(k^*)w$, and the fact that $A(k^*)>0$, which implies the least
eigenvalue $\lambda_1>0$. Based on the boundedness of the Rayleigh quotient of a
symmetric matrix, we have
\begin{align*}
\dot{V}=-w^T A(k^*)w\leq-\lambda_1\|w\|^2,
\end{align*}
which gives $\|w\|^2\leq-\dot{V}/\lambda_1$.  By integrating both sides, one has

\begin{align*}
\int_0^k\|w\|^2\,dt &\leq\int_0^k-\dot{V}/\lambda_1\,dt\\
&=V(0)/\lambda_1-V(k)/\lambda_1\\
&\leq V(0)/\lambda_1 \quad \text{for all } k. 
\end{align*}
This implies that the monotone increasing sequence
${\int_0^k\|w\|^2\,dt}$ is bounded above.  Hence,
$\int_0^\infty\|w\|^2\,dt$ converges.  According to Barb\u{a}lat's lemma,
$\lim_{t\to\infty}w(t)=0$. Therefore, the extended closed-loop system \eqref{eq:aug_system} is globally asymptotically stable. 
\end{proof}
\end{theorem}

\begin{figure}[H]
\centering
\subfloat[\centering]{\includegraphics[width=7.0cm]{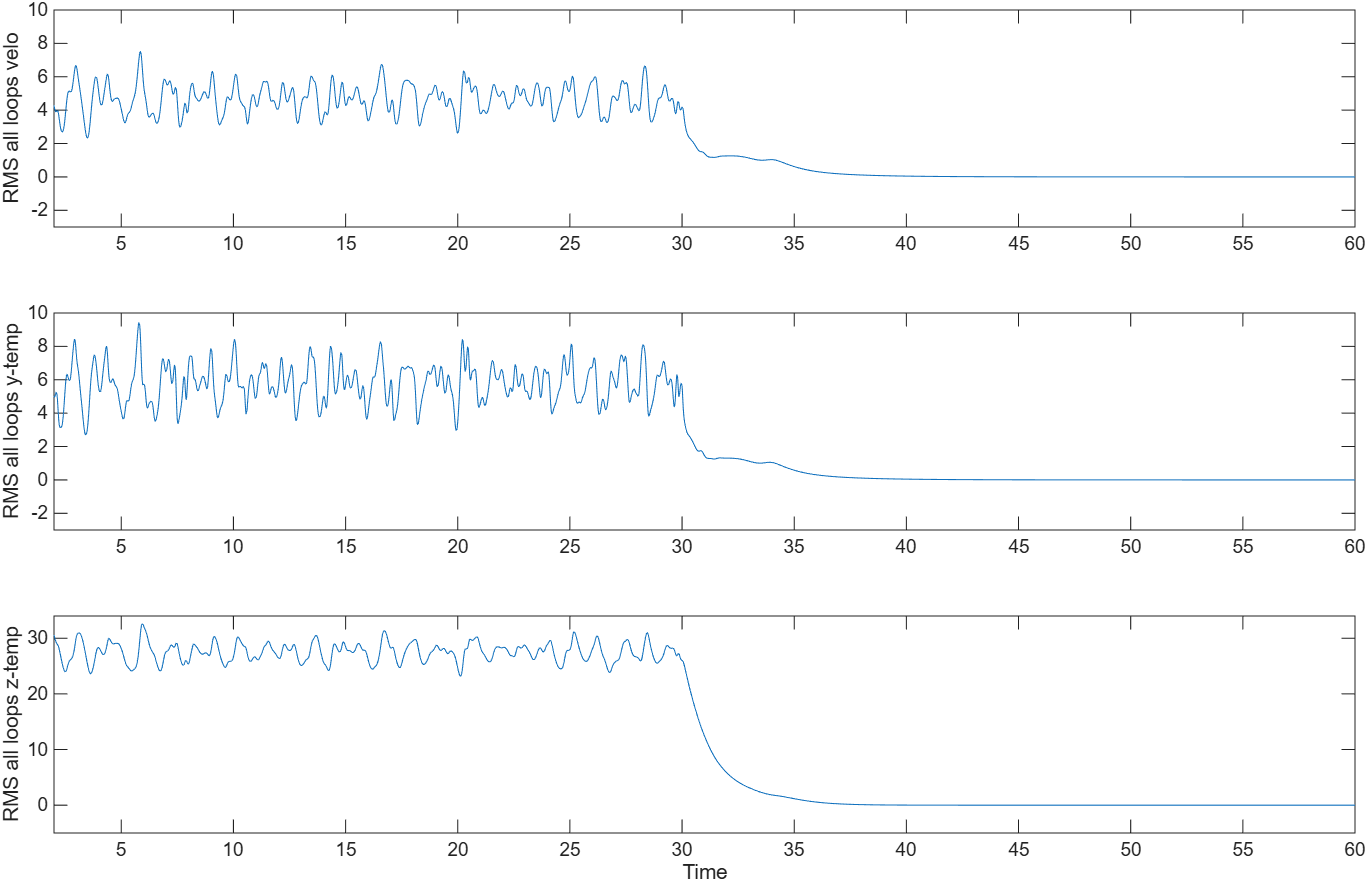}}
\subfloat[\centering]{\includegraphics[width=7.0cm]{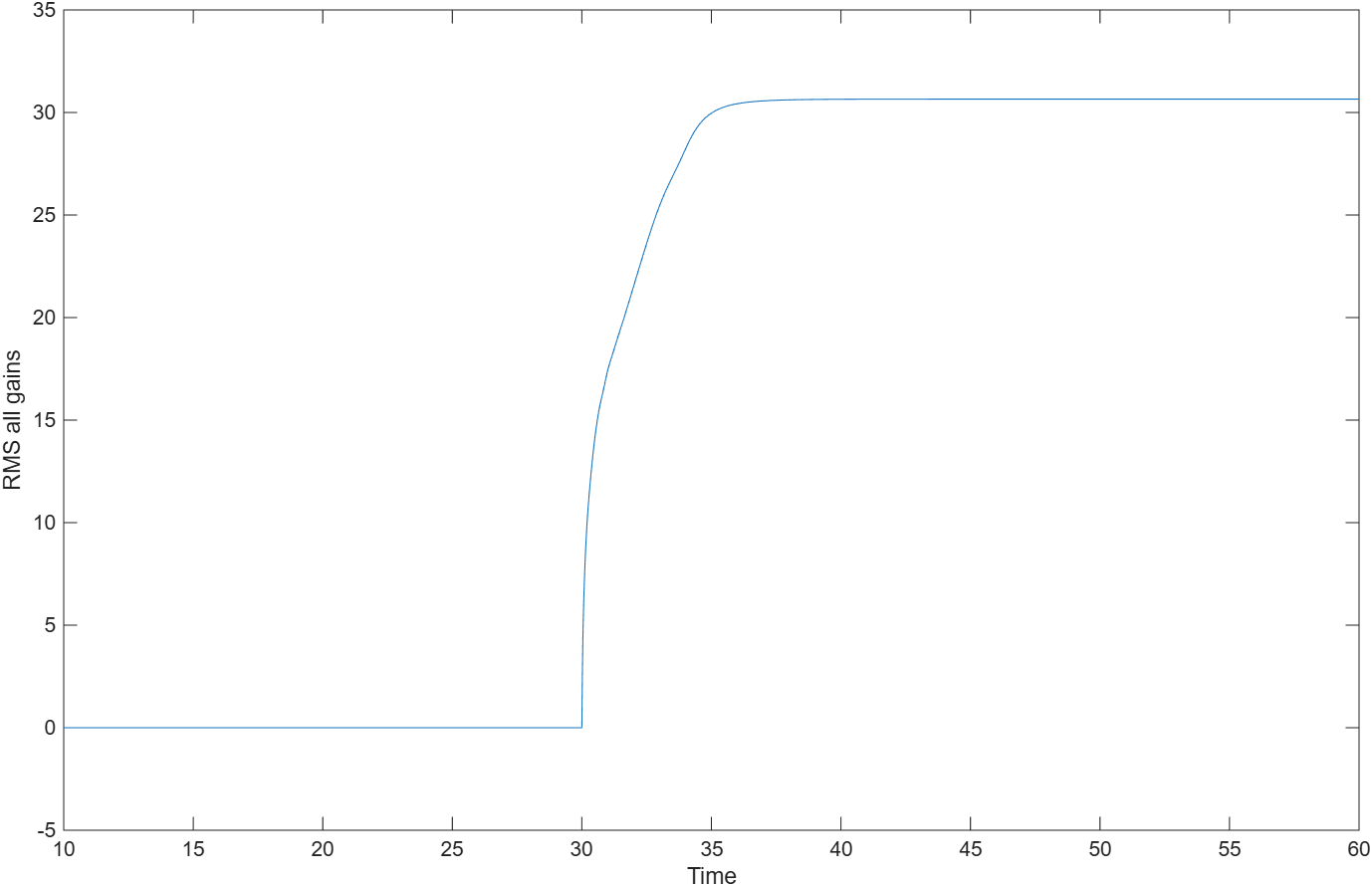}}
\caption{(\textbf{a}) Time response of quadratic mean of $x, y, z$-state with controller activated at $t=30s$ (\textbf{b}) Time response of quadratic mean of feedback gains with learning rate at 0.9 .\label{fig2}}
\end{figure}

While Theorem \ref{thm:positive_definiteness_before_recursions} and \ref{thm:global} showed that the decentralized adaptive control of the $N$-loop system \eqref{eq:aug_system} is globally asymptotically stable at the origin, the performance of the controller through the
transient and steady state responses are best observed from the numerical simulations. The
simulation results shown in Figure \ref{fig2} are based on a 30-loop system with 90 states and 90
differential equations. The range of the 30 Rayleigh numbers is between 25 and 55, the
momentum coupling parameter is set at 0.1, and the thermal coupling parameter is 0.2. Since
there are 90 states, it seems cumbersome to show the time response plots of even a subset of
those states as most of the state trajectories are similar in its chaotic nature. Instead, we choose
to show plots of the root-mean-square (RMS) or the quadratic means of the $x,y,z$-state, and adaptive gains $k_i$, respectively, at each sampling time $t$, as follows:
\begin{align}\label{eq:rms}
 x_r(t)=\sqrt{\frac{1}{N}\sum_{i=1}^{N}x_i^2(t)},\quad
& y_r(t)=\sqrt{\frac{1}{N}\sum_{i=1}^{N}y_i^2(t)},\quad
z_r(t)=\sqrt{\frac{1}{N}\sum_{i=1}^{N}z_i^2(t)},\notag\\ \quad
&k_r(t)=\sqrt{\frac{1}{N}\sum_{i=1}^{N}k_i^2(t)}
\end{align}
Such a measure by \eqref{eq:rms} is sufficient in demonstrating the time response of the states from
transient to steady state, particularly the steady state as the states converge to the origin. This is
readily seen from \eqref{eq:rms}  that each state of system \eqref{eq:aug_system}  is approaching 0 if the corresponding
quadratic mean \eqref{eq:rms} approaching 0, and vice versa. In all the simulations, it is assumed that
the system parameters are unknown. The controller gains evolve in time, which are determined
through the dynamic equations in \eqref{eq:aug_system}   with the initial value at zero. Figure \ref{fig2} (a) shows that the state trajectories converge to zero in the form of the quadratic mean of each of the $x$, $y$ and $z$-states. The controllers are activated at $t=30s$, and it takes about 8 s for the states to settle at the
equilibrium all together, with the learning rate set at $\beta_i=0.9$ in the dynamic equations of the gains in \eqref{eq:aug_system}.

Accordingly, the quadric mean of the gains approaches a steady state as shown in Figure \ref{fig2}(b). This implies that each gain $k_i$ in \eqref{eq:aug_system} converges to its own steady state $k_i^*$. We tested the matrix $A(k)$ evaluated at the steady state $k^*$, and the matrix $A(k^*)$ remains positive definite.

\subsection{Leading Principal Minors of the Reduced Block}
By Theorem~\ref{thm:positive_definiteness_before_recursions}, the
feedback gains can be chosen so that $H_N(k)\succ0
    \quad\text{and}\quad
    T_N(\eta)\succ0.$ Since $A_N(k) =
    \begin{bmatrix}
        H_N(k) & 0\\
        0 & T_N(\eta)
    \end{bmatrix},$ it follows that $A_N(k)\succ0.$

We now look into the leading principal minors of $H_N(k)$. For
$m=1,\ldots,2N$, let
\begin{align}\label{eq:leading_minor_definition}
\begin{cases}
    M_m
    &=
    H_N(k)[1:m,1:m],\\
    \mathcal A_m
    &=
    \det(M_m).
\end{cases}
\end{align}
We also use the convention $\mathcal A_0=1.$ Here, $M_m$ denotes the leading $m\times m$ principal submatrix of $H_N(k)$, whereas $\mathcal A_m$ denotes its determinant. Since $H_N(k)$ is symmetric and positive definite, Sylvester's
criterion gives
\begin{align}\label{eq:all_leading_minors_positive}
    \mathcal A_m>0,
    \qquad m=1,\ldots,2N.
\end{align}
The submatrix $H_N(k)$ is structured according to $(x_1,y_1,x_2,y_2,\ldots,x_N,y_N).$ For each $n=1,\ldots,N$, the minor $\mathcal A_{2n}$ corresponds to the states $(x_1,y_1,\ldots,x_n,y_n).$ For each $n=1,\ldots,N-1$, the minor $\mathcal A_{2n+1}$ corresponds to the states $(x_1,y_1,\ldots,x_n,y_n,x_{n+1}),$ whereas $\mathcal A_{2n+2}$ corresponds to the states $(x_1,y_1,\ldots,x_n,y_n,x_{n+1},y_{n+1}).$ Recall that $d_i$ is the diagonal entry of $H_N(k)$ corresponding to $x_i$, while $-q_i$ is the coupling entry between $x_i$ and $x_{i+1}$. 

The first two leading principal minors of $H_N(k)$ are
\begin{align}\label{eq:first_two_leading_minors}
\begin{cases}
    \mathcal A_1
    &=
    d_1
    =
    R_1(1+\gamma),\\
    \mathcal A_2
    &=
    \det
    \begin{bmatrix}
        d_1 & -R_1\\
        -R_1 & k_1
    \end{bmatrix}=
    d_1k_1-R_1^2.
\end{cases}    
\end{align}
The remaining principal minors can be generated recursively by
adding one state at a time: first $x_{n+1}$, and then
$y_{n+1}$. Theorem \ref{thm:positive_definiteness_before_recursions}
establishes that $H_N(k)\succ0$. Hence, all leading principal minors
of $H_N(k)$ are strictly positive. The following corollary records the
two recursive relations that will subsequently be used to extract
embedded lower-bound conditions for the feedback gains of an $N$-loop system.

\begin{corollary}\label{cor:recursive_relations_leading_minors}
The principal minors of $H_N(k)$ satisfy
\begin{align}\label{eq:recursion_A_2n_plus_1}
    \mathcal A_{2n+1}
    =
    d_{n+1}\mathcal A_{2n}
    -
    k_n\Delta_{n,n+1}\mathcal A_{2n-2},
    \qquad n=1,\ldots,N-1,
\end{align}
and
\begin{align}\label{eq:recursion_A_2n_plus_2_compact}
    \mathcal A_{2n+2}
    =
    k_{n+1}\mathcal A_{2n+1}
    -
    R_{n+1}^2\mathcal A_{2n},
    \qquad n=1,\ldots,N-1,
\end{align}
where
\begin{align}\label{eq:Delta_definition}
    \Delta_{n,n+1}
    =
    q_n^2
    = (R_n+R_{n+1})^2 \; \gamma^2/4,
    \qquad n=1,\ldots,N-1.
\end{align}

\begin{proof}
Fix $n\in\{1,\ldots,N-1\}$. The matrix $M_{2n+1}$ is obtained from
$M_{2n}$ by adjoining the variable $x_{n+1}$. Since $x_{n+1}$ couples
only with $x_n$ through the entry
$-q_n=-(R_n+R_{n+1})\gamma/2$, a cofactor expansion of $M_{2n+1}$ gives
\begin{align*}
    \det(M_{2n+1})
    =
    d_{n+1}\det(M_{2n})
    -
    q_n^2\det\left(M_{2n}^{(x_n)}\right),
\end{align*}
where $M_{2n}^{(x_n)}$ is obtained by deleting the row and column
corresponding to $x_n$. These deletions leave $y_n$ as an isolated
variable with diagonal entry $k_n$, while the remaining block is
$M_{2n-2}$. Hence,
\begin{align*}
    \det\left(M_{2n}^{(x_n)}\right)
    =
    k_n\det(M_{2n-2})
    =
    k_n\mathcal A_{2n-2}.
\end{align*}
Using $\det(M_{2n+1})=\mathcal A_{2n+1}$,\;
$\det(M_{2n})=\mathcal A_{2n}$, and
$q_n^2=\Delta_{n,n+1}$ by \eqref{eq:Delta_definition} to give
\begin{align*}
    \mathcal A_{2n+1}
    =
    d_{n+1}\mathcal A_{2n}
    -
    k_n\Delta_{n,n+1}\mathcal A_{2n-2}.
\end{align*}

Next, $M_{2n+2}$ is obtained from $M_{2n+1}$ by adjoining the
variable $y_{n+1}$. Since $y_{n+1}$ couples only with $x_{n+1}$
through the entry $-R_{n+1}$, a cofactor expansion of $M_{2n+2}$ gives
\begin{align*}
    \det(M_{2n+2})
    =
    k_{n+1}\det(M_{2n+1})
    -
    R_{n+1}^2
    \det\left(M_{2n+1}^{(x_{n+1})}\right).
\end{align*}
Deleting the row and column corresponding to $x_{n+1}$ leaves precisely $M_{2n}$. Therefore,
\begin{align*}
    \mathcal A_{2n+2}
    = k_{n+1}\mathcal A_{2n+1} -
    R_{n+1}^2\mathcal A_{2n}.\qedhere
\end{align*}
\end{proof}
\end{corollary}

The recursive relations in
Corollary \ref{cor:recursive_relations_leading_minors} can now be used
to extract explicit, embedded lower-bound expressions for the feedback gains. The iterations begin with $\mathcal A_0=1$, $\mathcal A_1=d_1$, and $\mathcal A_2=d_1k_1-R_1^2$. For $n=1$, one has $\mathcal A_2
    =
    k_1\mathcal A_1
    -
    R_1^2\mathcal A_0.$ For $n=2,\ldots,N$, replacing $n$ by $n-1$ in
\eqref{eq:recursion_A_2n_plus_2_compact} gives
\begin{align}\label{eq:general_even_minor_gain_identity}
    \mathcal A_{2n}
    =
    k_n\mathcal A_{2n-1}
    -
    R_n^2\mathcal A_{2n-2},
    \qquad n=1,\ldots,N.
\end{align}
For $n=1,\ldots,N-1$, substituting
\eqref{eq:general_even_minor_gain_identity} into
\eqref{eq:recursion_A_2n_plus_1} gives
\begin{align}\label{eq:odd_minor_factored_in_gain}
    \mathcal A_{2n+1}
    =
    k_n
    \left(
        d_{n+1}\mathcal A_{2n-1}
        -
        \Delta_{n,n+1}\mathcal A_{2n-2}
    \right)
    -
    d_{n+1}R_n^2\mathcal A_{2n-2}.
\end{align}
Next, we show that the coefficient of $k_n$ in
\eqref{eq:odd_minor_factored_in_gain} is strictly positive.
Under the hypotheses and feedback-gain conditions of
Theorem~\ref{thm:positive_definiteness_before_recursions},
we have $k_n>0$ and $H_N(k)\succ0$. Hence, Sylvester's
criterion gives
\begin{align*}
    \mathcal A_{2n+1}>0
    \qquad\text{and}\qquad
    \mathcal A_{2n-2}>0,
\end{align*}
where the second inequality also holds for $n=1$ by the
convention $\mathcal A_0=1$. Moreover, $d_{n+1}>0$ and
$R_n>0$. Consequently, both terms in $\mathcal A_{2n+1}
    +d_{n+1}R_n^2\mathcal A_{2n-2}$ are strictly positive. Adding $d_{n+1}R_n^2\mathcal A_{2n-2}$ to both sides of
\eqref{eq:odd_minor_factored_in_gain} and then dividing
by $k_n>0$, we obtain
\begin{align}\label{eq:positive_nested_bound_denominator}
    &d_{n+1}\mathcal A_{2n-1}
    -
    \Delta_{n,n+1}\mathcal A_{2n-2} =
    \frac{
        \mathcal A_{2n+1}
        +
        d_{n+1}R_n^2\mathcal A_{2n-2}
    }{k_n}
    >0.
\end{align}
Therefore, the coefficient of $k_n$ is strictly positive,
and division by this coefficient is valid when extracting
the lower bound on $k_n$. Hence, by setting \eqref{eq:odd_minor_factored_in_gain} greater than zero we obtain, the nonterminal feedback gains satisfy
\begin{align}\label{eq:nonterminal_recursive_gain_bound}
    k_n
    &>
    \frac{
        d_{n+1}R_n^2\mathcal A_{2n-2}
    }{
        d_{n+1}\mathcal A_{2n-1}
        -
        \Delta_{n,n+1}\mathcal A_{2n-2}
    }
    \notag\\
    &=
    \frac{
        R_n^2
    }{
        \displaystyle
        \frac{\mathcal A_{2n-1}}{\mathcal A_{2n-2}}
        -
        \frac{\Delta_{n,n+1}}{d_{n+1}}
    }
    >0,
    \qquad n=1,\ldots,N-1.
\end{align}
For the last gain $k_N$, \eqref{eq:general_even_minor_gain_identity}
and $\mathcal A_{2N}>0$ give
\begin{align}\label{eq:terminal_recursive_gain_bound}
    k_N
    >
    \frac{
        R_N^2\mathcal A_{2N-2}
    }{
        \mathcal A_{2N-1}
    }
    =
    \frac{
        R_N^2
    }{
        \displaystyle
        \frac{\mathcal A_{2N-1}}{\mathcal A_{2N-2}}
    }
    >0.
\end{align}
The ratios appearing in
\eqref{eq:nonterminal_recursive_gain_bound} and
\eqref{eq:terminal_recursive_gain_bound} can be expanded by repeated backward substitutions. Indeed, for $n=2,\ldots,N$, the two recursive relations give
\begin{align}\label{eq:minor_ratio_back_substitution}
    \frac{\mathcal A_{2n-1}}{\mathcal A_{2n-2}}
    =
    d_n
    -
    \frac{
        \Delta_{n-1,n}
    }{
        \displaystyle
        \frac{\mathcal A_{2n-3}}{\mathcal A_{2n-4}}
        -
        \frac{R_{n-1}^2}{k_{n-1}}
    },
\end{align}
with the initial ratio $\frac{\mathcal A_1}{\mathcal A_0}=d_1.$ Repeated substitutions of
\eqref{eq:minor_ratio_back_substitution} terminates at this initial
ratio and, therefore, produces a finite embedded expression for every feedback-gain bound.

\section{Nonlinear Thermal Coupling and Stability}
\label{sec:nonlinear_thermal_coupling}
It is concluded from the proof of
Theorem \ref{thm:positive_definiteness_before_recursions} that the
stability bounds on the feedback gains $k_i$ are independent of the
thermal coupling intensity parameter $\eta$ as long as $\eta>0$.
This has a lot to do with the $z$-states that are somehow decoupled
from the $x,y$-states by way of the stability matrix $A_N(k)$ associated with
$\dot V$. Because of this property, we investigated the possibility
of extending the result to more general settings such as a negative
lower bound on $\eta$ in the linear case, and the system with
nonlinear thermal coupling.

First off, we consider the scenario where $\eta$ is a constant. In this case, the matrix $T_N(\eta)$ is
written as
\begin{align}\label{eq:linear_thermal_matrix_decomposition}
    T_N(\eta)
    =
    I_N+\eta L_N,
\end{align}
where $L_N$ is given by
\begin{align}\label{eq:linear_thermal_laplacian}
    L_N
    =
    \begin{bmatrix}
        1  & -1 & 0  & 0 & \cdots & 0\\
        -1 & 2  & -1 & 0 & \cdots & 0\\
        0  & -1 & 2  & -1 & \cdots & 0\\
        \vdots & \vdots & \ddots & \ddots & \ddots & \vdots\\
        0 & 0 & \cdots & 0 & -1 & 1
    \end{bmatrix}_{N\times N}.
\end{align}

Our goal is to find a bound on $\eta$ so that $T_N(\eta)$ is
positive definite. This is achieved through the requirement that all
eigenvalues of $T_N(\eta)$ are greater than zero. Assume
$\lambda_j$ is an eigenvalue of $L_N$, then
$1+\eta\lambda_j$ is the corresponding eigenvalue of
$T_N(\eta)$. It is easy to see that $L_N$ is positive definite.
As such, $1+\eta\lambda_j>0$ yields a lower bound on $\eta$, i.e., $\eta>-\frac{1}{\lambda_j},
    \quad j=1,2,\ldots,N,$ which is guaranteed for all $j$ if $\eta>-\frac{1}{\lambda_{\max}},$ where $\lambda_{\max}$ is the greatest eigenvalue of $L_N$. To find the eigenvalues of $L_N$, we begin with
\begin{align*}
    L_Nv=\lambda v,
\end{align*}
where $v =
    [v_1\;v_2\;\cdots\;v_N]^\top.$ The eigenvalue problem yields the following system of equations for
$v$:
\begin{align}\label{eq:linear_thermal_eigenvalue_system}
    \begin{cases}
        v_1-v_2=\lambda v_1,\\
        -v_{j-1}+2v_j-v_{j+1}=\lambda v_j,
        \qquad j=2,3,\ldots,N-1,\\
        -v_{N-1}+v_N=\lambda v_N.
    \end{cases}
\end{align}
Let $v_j = \cos\left(j-\frac{1}{2}\right)\theta$,
After plugging this expression into
\eqref{eq:linear_thermal_eigenvalue_system}, one has
\begin{align*}
    -v_{j-1}+2v_j-v_{j+1}
    &=
    -\cos\left(j-\frac{3}{2}\right)\theta
    +
    2\cos\left(j-\frac{1}{2}\right)\theta-
    \cos\left(j+\frac{1}{2}\right)\theta\\
    &=
    4\sin^2\left(\frac{\theta}{2}\right)
    \cos\left(j-\frac{1}{2}\right)\theta.
\end{align*}
Therefore, $\lambda
    = 4\sin^2\left(\frac{\theta}{2}\right).$ To find the distinct eigenvalues of $L_N$, we first plug the
expressions into the first equation in
\eqref{eq:linear_thermal_eigenvalue_system}, which only yields an
identity. However, the last equation gives
\begin{align*}
    &-\cos\left(N-\frac{3}{2}\right)\theta
    +
    \cos\left(N-\frac{1}{2}\right)\theta\\
    &\qquad=
    4\sin^2\left(\frac{\theta}{2}\right)
    \cos\left(N-\frac{1}{2}\right)\theta,
\end{align*}
which is simplified to be $\sin\left(\frac{\theta}{2}\right)\sin(N\theta)=0.$ Obviously, $\sin\left(\frac{\theta}{2}\right)\neq0,$ otherwise, all the eigenvalues of $L_N$ would be $0$, which is
false. From $\sin(N\theta)=0,$ we have $\theta_j = \frac{(j-1)\pi}{N}, \quad 1\leq j\leq N.$
Hence, the greatest eigenvalue of $L_N$ is
\begin{align*}
    \lambda_{\max}
    &=
    4\sin^2\left(
        \frac{(N-1)\pi}{2N}
    \right)\\
    &=
    2+2\cos\left(\frac{\pi}{N}\right).
\end{align*}
We then arrive at the conclusion that $T_N(\eta)$ is positive
definite if
\begin{align}\label{eq:constant_eta_lower_bound}
    \eta>
    -\frac{1}{2+2\cos(\pi/N)},
\end{align}
which serves as the lower bound on the coupling intensity constant
between the $z$-states. This constraint is automatically satisfied
if $\eta>0$.

Next, we consider the case of nonlinear thermal coupling by way of
replacing the thermal coupling terms by nonlinear functions that
depend on the temperature difference between the adjacent loops. The
modified $z$-equations are given as follows:
\begin{align}\label{eq:nonlinear_z_equations_f}
    \begin{cases}
        \dot z_1
        =
        x_1y_1-z_1-f(z_1-z_2),\\
        \dot z_2
        =
        x_2y_2-z_2-f(z_2-z_1)-f(z_2-z_3),\\
        \hspace{8mm}\vdots\\
        \dot z_{N-1}
        =
        x_{N-1}y_{N-1}-z_{N-1}
        -f(z_{N-1}-z_{N-2})
        -f(z_{N-1}-z_N),\\
        \dot z_N
        =
        x_Ny_N-z_N-f(z_N-z_{N-1}).
    \end{cases}
\end{align}
Here, $f:S\rightarrow\mathbb R,$ where $S$ is a compact subset of $\mathbb R$ because the state
trajectory stays on its bounded attractor. It is reasonable to assume that
\begin{align*}
    f(0)=0
\end{align*}
since $f$ prescribes the heat transferred between the adjacent loops due to their temperature difference at the coupling point. Moreover,
according to the second law of thermodynamics, heat transfers from
warmer body to cooler body. As such, $f(x)$ can be written as
\begin{align}\label{eq:f_factorization_original}
    f(x)=xg(x),
\end{align}
where $g(-x)=g(x).$ Now, system \eqref{eq:nonlinear_z_equations_f} has the form
\begin{align}\label{eq:nonlinear_z_equations_g}
    \begin{cases}
        \dot z_1
        =
        x_1y_1-z_1-(z_1-z_2)g_1,\\
        \dot z_2
        =
        x_2y_2-z_2-(z_2-z_1)g_1-(z_2-z_3)g_2,\\
        \hspace{8mm}\vdots\\
        \dot z_{N-1}
        =
        x_{N-1}y_{N-1}-z_{N-1}
        -(z_{N-1}-z_{N-2})g_{N-2}
        -(z_{N-1}-z_N)g_{N-1},\\
        \dot z_N
        =
        x_Ny_N-z_N-(z_N-z_{N-1})g_{N-1},
    \end{cases}
\end{align}
where
\begin{align}\label{eq:gi_original_definition}
    g_i
    =
    g(z_i-z_{i+1})
    =
    g(z_{i+1}-z_i),
    \qquad i=1,2,\ldots,N-1.
\end{align}
After substituting the $z$-equations with
\eqref{eq:nonlinear_z_equations_g}, we study the stability of the
following coupled $N$-loop system with nonlinear thermal coupling:
\begin{align}\label{eq:nonlinear_thermal_full_system_original}
    \begin{cases}
        \displaystyle
        \dot x_i
        =
        p\left[
            y_i-x_i
            -\gamma_{i-1}(x_i-x_{i-1})
            -\gamma_i(x_i-x_{i+1})
        \right],\\[2mm]
        \displaystyle
        \dot y_i
        =
        R_ix_i-x_iz_i-k_iy_i,\\
        \displaystyle
        \dot z_i
        =
        x_iy_i-z_i
        -(z_i-z_{i-1})g(z_i-z_{i-1})
        -(z_i-z_{i+1})g(z_i-z_{i+1}),
        \quad i=1,2,\ldots,N,
    \end{cases}
\end{align}
with the assumption $g(z_1-z_0)
    =  g(z_N-z_{N+1}) = 0$ for the redundant terms.

To perform the stability analysis of
\eqref{eq:nonlinear_thermal_full_system_original}, we consider the same Lyapunov function as before:
\begin{align}\label{eq:nonlinear_original_lyapunov_function}
    V
    =
    \frac{1}{2}
    \sum_{i=1}^{N}
    \left(
        \frac{R_i}{p}x_i^2+y_i^2+z_i^2
    \right).
\end{align}
Let $X
    =
    [x_1,y_1,x_2,y_2,\ldots,x_N,y_N,
    z_1,\ldots,z_N]^\top.$ The time derivative of $V$ along the state trajectory is written as
\begin{align}\label{eq:nonlinear_original_matrix_derivative}
    \dot V
    =
    -X^\top A_N^{\mathrm{nl}}(k,z)X,
\end{align}
where
\begin{align}\label{eq:nonlinear_original_block_matrix}
    A_N^{\mathrm{nl}}(k,z)
    =
    \begin{bmatrix}
        H_N(k) & 0\\
        0 & T_N^{\mathrm{nl}}(z)
    \end{bmatrix}.
\end{align}
The submatrix $H_N(k)$ is the same as the corresponding submatrix in
the linear thermal-coupling case, while the submatrix
$T_N^{\mathrm{nl}}(z)$ is given by
\begin{align}\label{eq:nonlinear_original_thermal_matrix}
    T_N^{\mathrm{nl}}(z)
    =
    \begin{bmatrix}
        1+g_1
        & -g_1
        & 0
        & \cdots
        & 0\\
        -g_1
        & 1+g_1+g_2
        & -g_2
        & \ddots
        & \vdots\\
        0
        & -g_2
        & 1+g_2+g_3
        & \ddots
        & 0\\
        \vdots
        & \ddots
        & \ddots
        & \ddots
        & -g_{N-1}\\
        0
        & \cdots
        & 0
        & -g_{N-1}
        & 1+g_{N-1}
    \end{bmatrix}.
\end{align}

According to
Theorem~\ref{thm:positive_definiteness_before_recursions}, the
submatrix $H_N(k)$ is positive definite if
\eqref{eq:theorem_rayleigh_closeness} and
\eqref{eq:theorem_preliminary_gain_bounds} are satisfied. Therefore,
the matrix $A_N^{\mathrm{nl}}(k,z)$ is positive definite if the
submatrix $T_N^{\mathrm{nl}}(z)$ is positive definite. We will find
the condition on $g(x)$ so that
\begin{align*}
    T_N^{\mathrm{nl}}(z)\succ0.
\end{align*}
To this end, given the structure of
\eqref{eq:nonlinear_original_thermal_matrix}, the matrix can be written as
\begin{align}\label{eq:nonlinear_original_factorization}
    T_N^{\mathrm{nl}}(z)
    =
    I_N+B^\top DB,
\end{align}
where $D  = \operatorname{diag}
    (g_1,g_2,\ldots,g_{N-1}),$
and the matrix $B$ is an $(N-1)\times N$ first-order difference
matrix given by
\begin{align}\label{eq:original_difference_matrix}
    B
    =
    \begin{bmatrix}
        1 & -1 & 0 & \cdots & 0\\
        0 & 1 & -1 & \ddots & \vdots\\
        0 & 0 & 1 & \ddots & 0\\
        \vdots & \ddots & \ddots & \ddots & -1\\
        0 & \cdots & 0 & 1 & -1
    \end{bmatrix}.
\end{align}
Let $v\in\mathbb R^N$, then
\begin{align}\label{eq:original_nonlinear_quadratic_form}
    v^\top T_N^{\mathrm{nl}}(z)v
    =
    \|v\|^2
    +
    \sum_{j=1}^{N-1}
    g_j(v_j-v_{j+1})^2.
\end{align}
Let $m  =    \min_{x\in S}g(x).$
If $m\geq0$, it is easy to see that $v^\top T_N^{\mathrm{nl}}(z)v>0$ because $g_j\geq0$. Therefore, we assume $m<0$. In this case,  $D\succeq mI_{N-1},$ from which one has
\begin{align*}
    v^\top T_N^{\mathrm{nl}}(z)v
    &=
    \|v\|^2+(Bv)^\top D(Bv)\\
    &\geq
    \|v\|^2+m(v^\top B^\top Bv).
\end{align*}
The matrix $B^\top B$ turns out to be $L_N$ given by
\eqref{eq:linear_thermal_laplacian}. Since $m<0$,
\begin{align*}
    m(v^\top B^\top Bv)
    \geq
    m\lambda_{\max}(B^\top B)\|v\|^2,
\end{align*}
where $\lambda_{\max}(B^\top B)$ is the greatest eigenvalue of
$B^\top B$, and
\begin{align*}
    \lambda_{\max}(B^\top B)
    =
    2+2\cos\left(\frac{\pi}{N}\right).
\end{align*}
Hence, we have
\begin{align*}
    v^\top T_N^{\mathrm{nl}}(z)v
    \geq
    \left[
        1+m
        \left(
            2+2\cos\left(\frac{\pi}{N}\right)
        \right)
    \right]\|v\|^2.
\end{align*}
Therefore, $v^\top T_N^{\mathrm{nl}}(z)v>0$ if
\begin{align*}
    1+m
    \left[
        2+2\cos\left(\frac{\pi}{N}\right)
    \right]
    >0,
\end{align*}
or
\begin{align*}
    m =\min_{x\in S}g(x) >
    -\frac{1}{2+2\cos(\pi/N)}.
\end{align*}
This implies a sufficient condition for
$T_N^{\mathrm{nl}}(z)$ to be positive definite, and that is
\begin{align}\label{eq:original_nonlinear_g_condition}
    g(x)
    >
    -\frac{1}{2+2\cos(\pi/N)},
    \qquad
    \text{for all }x\in S.
\end{align}

It is worth noting that the lower bound on the constant coupling
intensity parameter $\eta$ also serves as the baseline for the
nonlinear function $g$, and such a lower bound has an explicit form.

Finally, according to
Theorem~\ref{thm:positive_definiteness_before_recursions}, the
$N$-loop system
\eqref{eq:nonlinear_thermal_full_system_original} with nonlinear
thermal coupling is globally asymptotically stable if the feedback
gains $k_i$ satisfy
\eqref{eq:theorem_preliminary_gain_bounds} and the coupling function
$g$ satisfies \eqref{eq:original_nonlinear_g_condition}.

\begin{figure}[H]
\centering
\subfloat[\centering]{\includegraphics[width=7.0cm]{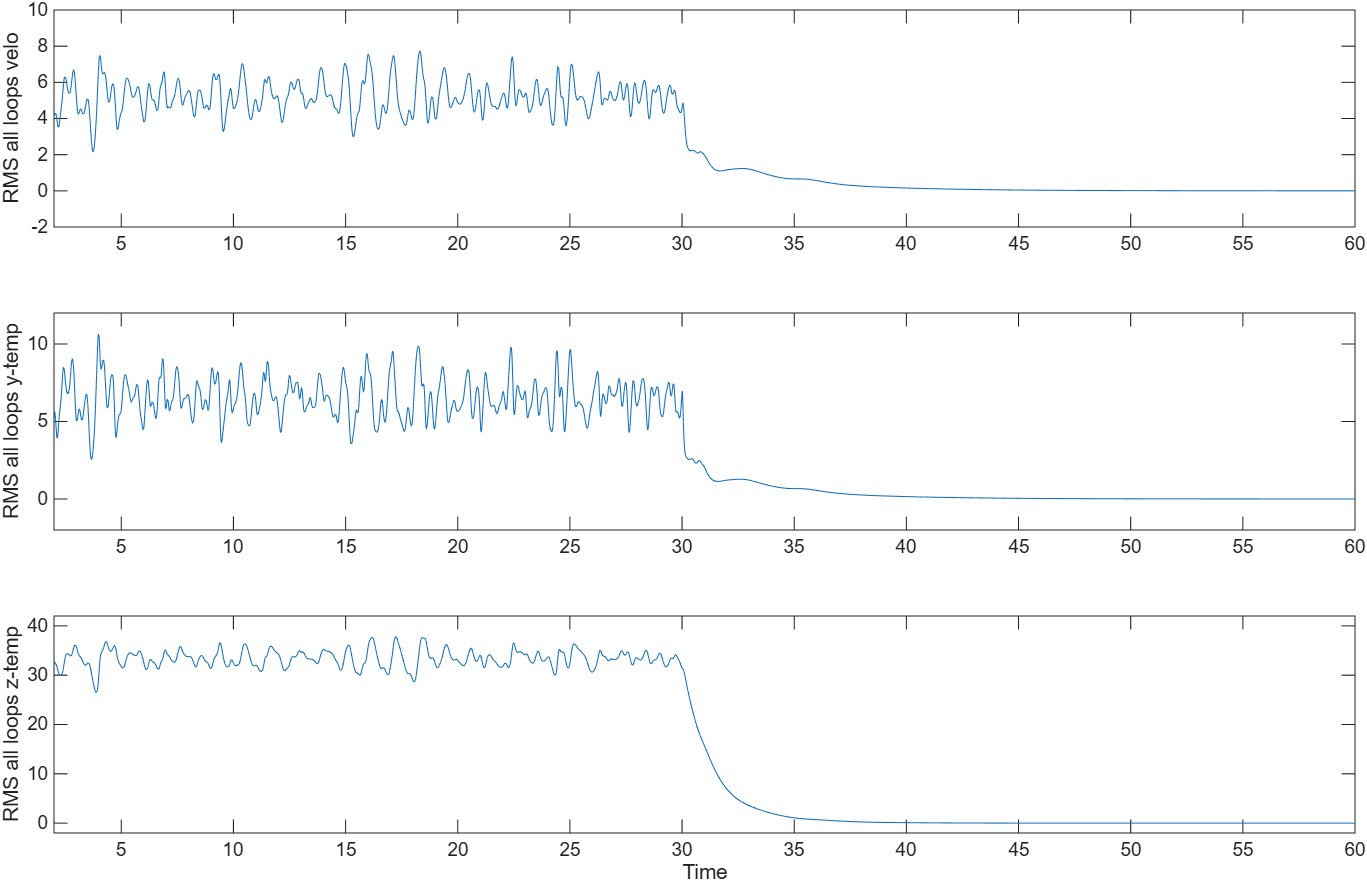}}
\subfloat[\centering]{\includegraphics[width=7.0cm]{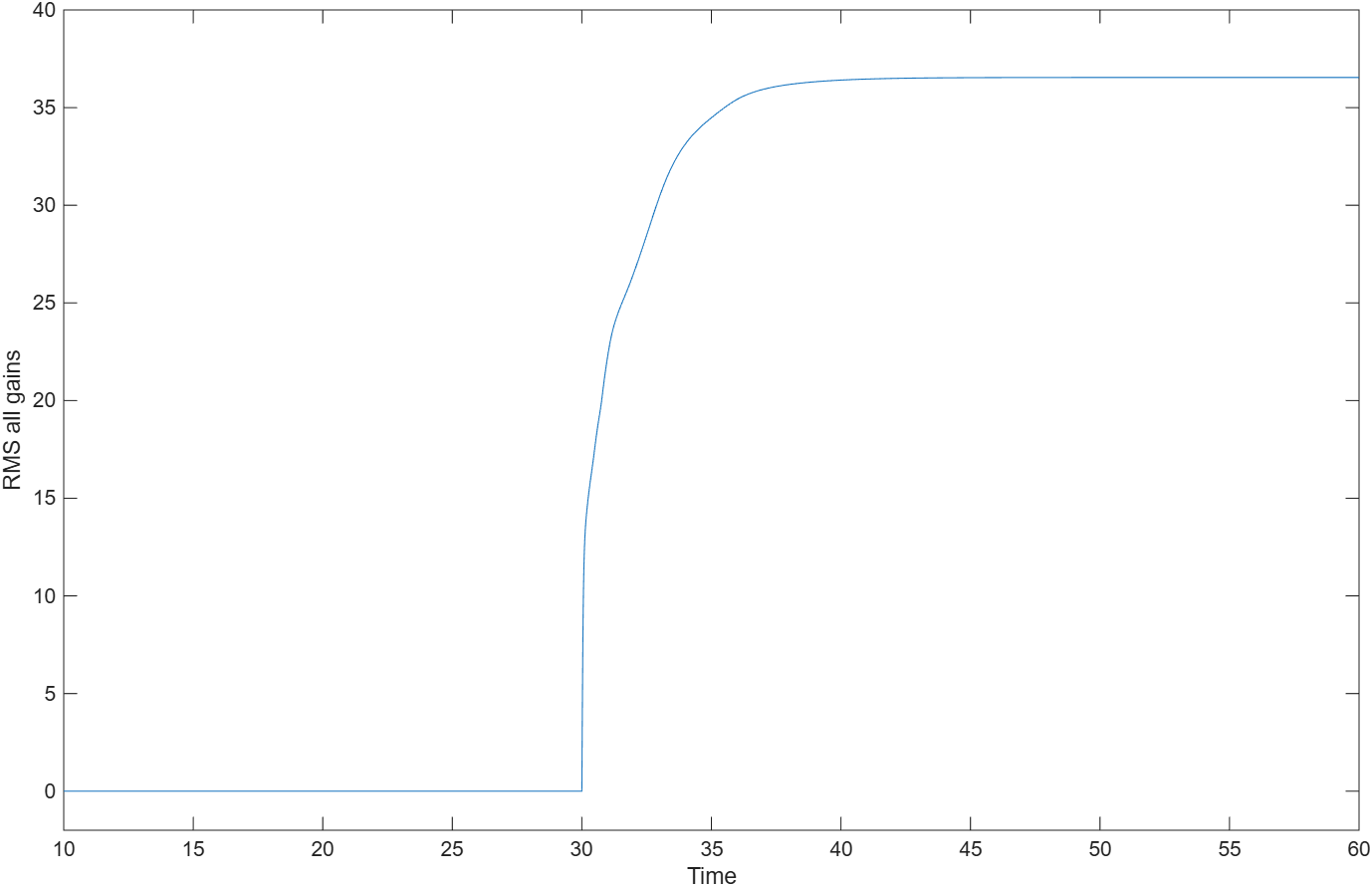}}
\caption{(\textbf{a}) Time response of quadratic mean of $x, y, z$-state with nonlinear thermal coupling (\textbf{b})Time response of quadratic mean of feedback gains under nonlinear thermal coupling.\label{fig3}}
\end{figure}

As discussed above, the thermal coupling has no impact on the stability of the control
system \eqref{eq:n_sys} if the coupling satisfies the conditions given by \eqref{eq:constant_eta_lower_bound}. This allows an extension to
nonlinear thermal coupling when the condition \eqref{eq:original_nonlinear_g_condition} is met.
We tested the system with the nonlinear function $f(z)=z(1+\cos z)$ that prescribes the thermal
coupling, which indeed satisfies the constraint \eqref{eq:original_nonlinear_g_condition}. Using the same values for the Rayleigh
numbers, momentum coupling constants, and the learning rate of the adaptive gains as with the
previous simulations, we show in Figure \ref{fig3} that the simulation results agree with the theoretical
finding. In particular, the nonlinear coupling does impact the trajectories of the z-states as well as
the $x$, $y$-states. However, it does not influence the stability of system \eqref{eq:nonlinear_thermal_full_system_original}. Moreover, the transient
responses and settling time between Figures \ref{fig2}(a) and \ref{fig3}(a) are also comparable. Meanwhile, the
corresponding time response of the quadratic mean of the adaptive gains is shown in Figure \ref{fig3}(b).

\begin{figure}[H]{\centering}{} 
\includegraphics[width=9.0cm]{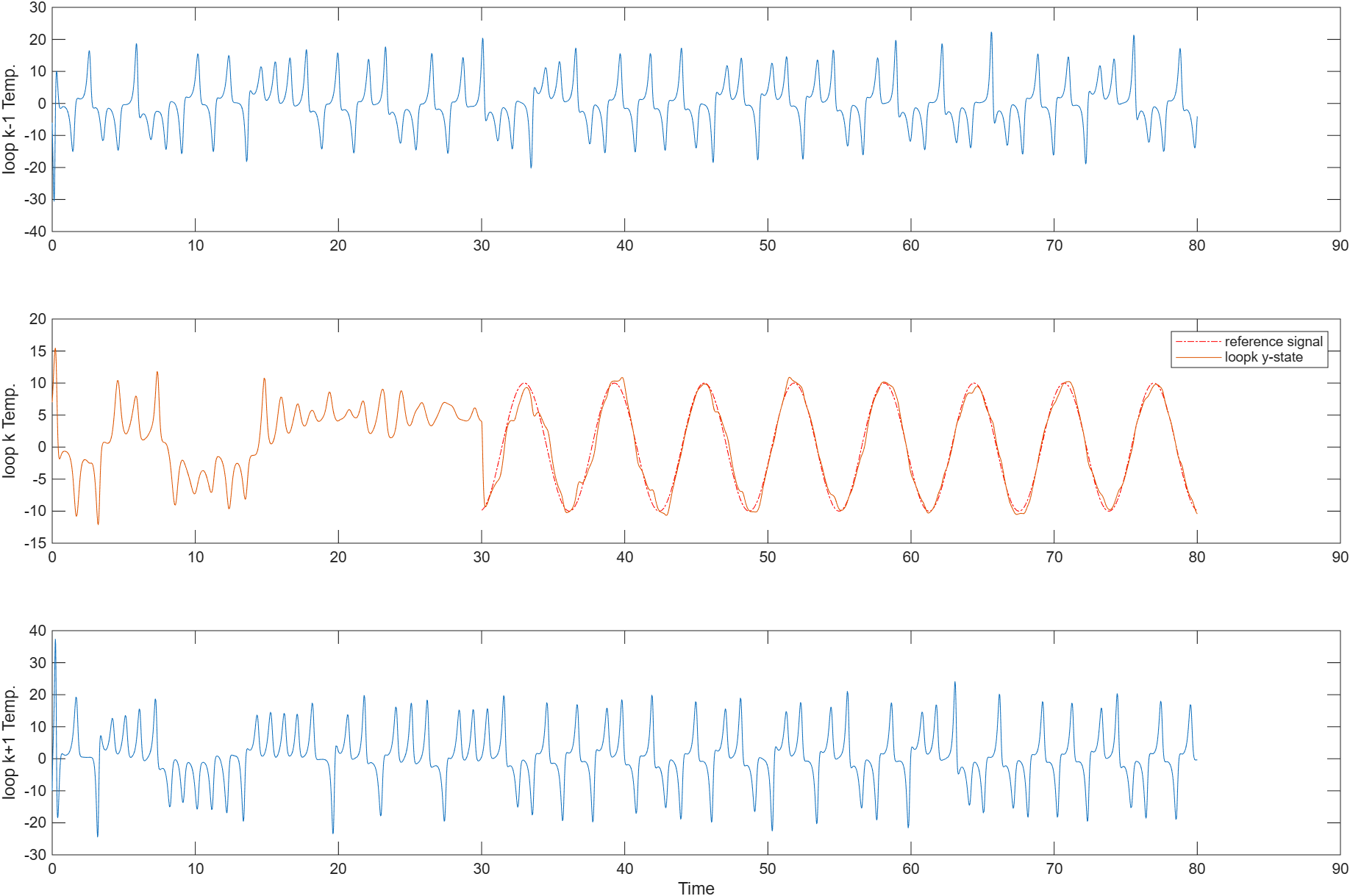}\caption{Time response of 3-loop system with middle loop tracking a sinusoid and chaotic
adjacent loops.}\label{fig4}
\end{figure}

\begin{figure}[H]{\centering}{} 
\includegraphics[width=9.0cm]{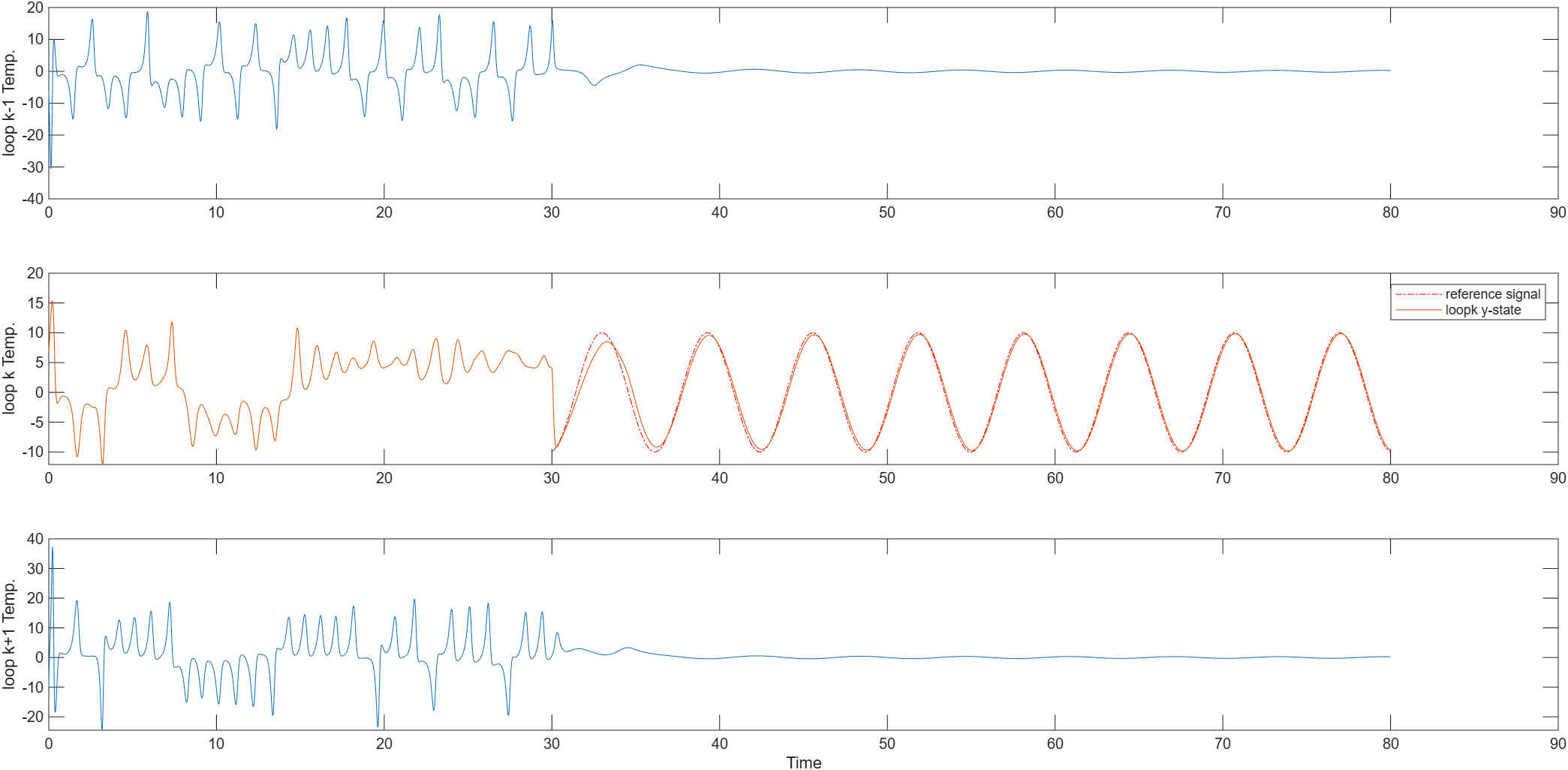}\caption{Time response of 3-loop system with middle loop tracking a sinusoid and stabilized
adjacent loops.}\label{fig5}
\end{figure}

Theorems \ref{thm:positive_definiteness_before_recursions} and \ref{thm:global} guarantee that the system \eqref{eq:n_sys} is globally asymptotically stable at the origin
so as shown from the simulations. We also tested the tracking ability of the system when the reference signal is a nonstationary signal, such as a sinusoid. To this end, we modify the adaptive control system \eqref{eq:aug_system} to take the reference signal into account as follows
\begin{align}\label{eq:track}
\begin{cases}
    \dot{x}_i &= p\{y_i-x_i-\gamma_{i-1}(x_i-x_{i-1})
    -\gamma_i(x_i-x_{i+1})\} \\
    \dot{y}_i &= R_i x_i-y_i-x_i z_i-k_i(y_i-r_i) \\
    \dot{z}_i &= x_i y_i-z_i-\eta(z_i-z_{i-1})
    -\eta(z_i-z_{i+1}) \\
    \dot{k}_i &= \beta_i(y_i-r_i)^2,\ i=1,2,\ldots,n
    \end{cases}
\end{align}
where $r_i$ is the reference signal. One of the interesting cases for testing the tracking performance of \eqref{eq:track}
is to make one of the middle loops track an input signal while other loops remain chaotic. We demonstrate the tracking ability via a three-loop system. Again, we use large Rayleigh numbers to
drive the system into chaos and then activate tracking control only for the middle loop. Figure \ref{fig4} shows the time response of the $y$-state of each loop. In this experiment, the reference signal is a sine wave. Clearly, the $y$-state of the middle loop is tracking the sinusoid while its adjacent loops
are chaotic. The tracking error is small, but best described as bounded-input, bounded-output
stable tracking with a steady state error. Such an error is caused by continuous interference from the neighboring loops by way of the momentum coupling and thermal coupling. On the other hand, the tracking accuracy is significantly improved if the state trajectories of the adjacent loops are stabilized at the equilibrium, see Figure \ref{fig5}.

\section{Conclusions}\label{sec:conclusions}
In this paper, we developed a decentralized stabilization
framework for a laterally coupled thermosyphon system with
an arbitrary finite number of loops. By separating the
thermal states from the remaining states, the stability
matrix is decoupled into
in block-diagonal form. This structure allows us to establish the existence of lower bounds on the feedback
gains and prove the global asymptotic stability of the
closed-loop system. This leads to the adaptive controller design when the system parameters are unknown. The two recursive relations for the
leading principal minors provide a systematic way to
extract explicit gain constraints and computable bounds on the gains. The block decomposition of the stability matrix also reveals the role of thermal
coupling in the stability analysis as the stability bounds on the feedback gains do not depend on the thermal coupling intensity constant. The thermal-block submatrix
gives an explicit admissible range for
negative linear coupling and a sufficient condition for
a class of nonlinear thermal coupling functions. These
extensions preserve stability, although thermal coupling
can still affect the transient response.

\bibliographystyle{unsrtnat}
\bibliography{sample}

@article{ref1,
  author  = {Kordestani, M. and Safavi, A. A. and Saif, M.},
  title   = {{Recent Survey of Large-Scale Systems: Architectures, Controller Strategies and Industrial Applications}},
  journal = {IEEE Syst. J.},
  year    = {2021},
  volume  = {15},
  pages   = {5440--5453}
}

@article{ref2,
  author  = {Li, Y. and Xu, B. and Ahn, C.},
  title   = {{A Cyclic-Small-Gain-Based Adaptive Fuzzy Control of Interconnected Systems}},
  journal = {IEEE Syst. J.},
  year    = {2022},
  volume  = {16},
  pages   = {6468--6479}
}

@article{ref3,
  author  = {Ehrhard, P. and Karcher, C. and Muller, U.},
  title   = {{Dynamical Behavior of Natural Convection in a Double-Loop System}},
  journal = {Expri. Heat Transf.},
  year    = {1989},
  volume  = {2},
  pages   = {13--26}
}

@article{ref4,
  author  = {Wu, Y. and Braselton, J. and Jin, Y. and Shahat, A.},
  title   = {{Adaptive Control of Bi-Directionally Coupled Lorenz Systems with Uncertainties}},
  journal = {J. Frankl. Inst.},
  year    = {2019},
  volume  = {356},
  pages   = {1287--1301}
}

@misc{ref5,
  author = {Anderson, N.},
  title  = {{Stabilize Chaotic Flows in a Coupled Triple-Loop Thermosyphon System}},
  year   = {2019},
  note   = {Honors College Thesis, 452}
}

@misc{ref6,
  author = {Dey, N. K.},
  title  = {{Stability Control of Laterally Coupled Quad-Loop Dynamical Systems}},
  year   = {n.d.},
  note   = {College of Graduate Thesis \& Dissertation, 2532}
}

@article{ref7,
  author  = {Dey, N. K. and Wu, Y.},
  title   = {{Adaptive Control of Laterally Coupled Quadruple Lorenz Systems}},
  journal = {Nonlinear Stud.},
  year    = {2024},
  volume  = {31},
  pages   = {1037--1051}
}

@article{ref8,
  author  = {Dey, N. K. and Wu, Y.},
  title   = {{Decentralized Disturbance Rejection Control of Triangularly Coupled Loop Thermosyphon System}},
  journal = {Actuators},
  year    = {2025},
  volume  = {14},
  pages   = {532}
}

@article{ref9,
  author  = {Sun, L. and Huang, X. and Song, Y.},
  title   = {{Decentralized Intermittent Feedback Adaptive Control of Non-Triangular Nonlinear Time-Varying Systems}},
  journal = {IEEE Trans. Autom. Control},
  year    = {2024},
  volume  = {69},
  pages   = {1265--1272}
}

@article{ref10,
  author  = {Yang, W. and Xia, J. and Yu, M. and Zhang, N.},
  title   = {{Decentralized Adaptive Funnel Control of Uncertain Large-Scale Interconnected Nonlinear System}},
  journal = {Appl. Math. Comput.},
  year    = {2023},
  volume  = {441},
  pages   = {127694}
}

@article{ref11,
  author  = {Jiang, Z. and Zhang, H. and Xue, L.},
  title   = {{A Semi-Global Finite-Time Decentralized Control Method for High-Order Large-Scale Nonlinear Systems}},
  journal = {Actuators},
  year    = {2024},
  volume  = {13},
  pages   = {250}
}

@article{ref12,
  author  = {He, W. and Liu, Y. and Zhang, Q.},
  title   = {{Decentralized Output-Feedback Adaptive Event-Triggered Control for Interconnected Nonlinear Delay Systems with Actuator Failures}},
  journal = {Actuators},
  year    = {2024},
  volume  = {13},
  pages   = {188}
}

@article{ref13,
  author  = {Li, M. and Ahn, C. K. and Xiang, Z.},
  title   = {{Decentralized Adaptive Fuzzy Finite-Time Event-Triggered Control for Interconnected Nonlinear Systems Subject to Input Saturation}},
  journal = {IEEE Syst. J.},
  year    = {2023},
  volume  = {17},
  pages   = {1648--1659}
}

@article{ref14,
  author  = {Machado, J. E. and Cucuzzella, M. and Pronk, N. and Scherpen, J. M. A.},
  title   = {{Adaptive Control for Flow and Volume Regulation in Multi-Producer District Heating Systems}},
  journal = {IEEE Control Syst. Lett.},
  year    = {2022},
  volume  = {6},
  pages   = {794--799}
}

@article{ref15,
  author  = {Machado, J. E. and Ferguson, J. and Cucuzzella, M. and Scherpen, J. M. A.},
  title   = {{Decentralized Temperature and Storage Volume Control in Multi-Producer District Heating}},
  journal = {IEEE Control Syst. Lett.},
  year    = {2023},
  volume  = {7},
  pages   = {413--418}
}

\end{document}